\documentstyle[12pt]{article}
\newlength{\abstractwidth}
\renewcommand{\thefootnote}{\fnsymbol{footnote}}
\renewcommand{\thanks}[1]{\footnote{#1}} 
\newcommand{\starttext}{
\setcounter{footnote}{0}
\renewcommand{\thefootnote}{\arabic{footnote}}}

\newcommand{\be}{\begin{equation}}
\newcommand{\bea}{\begin{eqnarray}}
\newcommand{\eea}{\end{eqnarray}}
\newcommand{\ee}{\end{equation}}

\def\ba{\begin{eqnarray}}
\def\ea{\end{eqnarray}}

\def\p{\partial}

\def\o{\omega}

\def\al{\alpha}

\def\e{\epsilon}

\def\o{\omega}

\def\s{\sigma}

\def\v{\vskip .1in}

\def\[{{\bf [}}
\def\]{{\bf ]}}

\begin{document}
\starttext
\baselineskip=18pt
\setcounter{footnote}{0}

\begin{center}
{\Large \bf ON ASYMPTOTICS FOR THE MABUCHI } \\
\bigskip
{\Large \bf ENERGY FUNCTIONAL}
\footnote{Research supported in part by National Science Foundation
grants DMS-98-00783 and DMS-01-00410}

\bigskip\bigskip

{\large D.H. Phong$^*$ and Jacob Sturm$^{\dagger}$} \\ 

\bigskip

$^*$ Department of Mathematics\\
Columbia University, New York, NY 10027\\
\v
$^{\dagger}$ Department of Mathematics \\
Rutgers University, Newark, NJ 07102

\end{center}

\baselineskip=15pt
\setcounter{equation}{0}
\setcounter{footnote}{0}

\section{Introduction}
\setcounter{equation}{0}

Many canonical geometric structures have been found
to be closely related to stability in the sense of
geometric invariant theory. For the existence of K\"ahler-Einstein metrics, this
is the classical conjecture of Yau \cite{yau, yau1, yau2}. The necessity of several notions of
stability has been established in this case by Tian \cite{tian97} and by Donaldson
\cite{donaldson}.
In the variational approach, the existence of K\"ahler-Einstein metrics can be
reduced to the behavior of energy functionals \cite{tian97}. Of particular importance
is the Futaki energy functional $F_{\o_0}^0(\phi)$ \cite{futaki,
ding}
\be
F_{\o_0}^0(\phi)
=-{1\over (n+1)V}\int_X\phi\sum_{i=0}^n\o_0^i\o_{\phi}^{n-i}
\ee
and the Mabuchi energy
functional $\nu_{\o_0}(\phi)$ \cite{mabuchi, chen}
\be
\label{mabuchi}
\nu_{\o_0}(\phi)
=
{1\over V}
\int_X\big\{(\ln{\o_{\phi}\over\o_0})\o_\phi^n
-
\phi(Ric(\o_0)\sum_{i=0}^{n-1}\o_0^i\o_\phi^{n-1-i}
-{\mu(X)\over n+1}\sum_{i=0}^n\o_0^i\o_\phi^{n-i})\big\}
\ee
Here $X$ is a compact complex manifold of dimension $n$,
$\o_0$ is a reference K\"ahler form, $\o_\phi=\o_0+{\sqrt{-1}\over 2\pi}
\p\bar\p\phi$, $Ric(\o_0)=-{\sqrt{-1}\over 2\pi}\p\bar\p
\,\ln\,\o_0^n$ is the
Ricci curvature of the K\"ahler form $\o_0$, $V$ is the volume
of $X$ with respect to $\o_0$, and $\mu(X)$ is the average
scalar curvature
\be
\mu(X)={1\over V}\int_Xn\,Ric(\o_0)\wedge \o_0^{n-1}
\ee
In terms of
the functional $J_{\o_0}(\phi)$ of Aubin and Yau,
$F_{\o_0}^0(\phi)$ can also be recast as
\be
\label{aubinyau}
F_{\o_0}^0(\phi)
=
{\sqrt{-1}\over 2\pi V}\int_X\sum_{i=0}^{n-1}
{(i+1)
\over (n+1)}
\p\phi\wedge\bar\p\phi\wedge
\o_{\phi}^{n-i-1}
\wedge
\o_0^i
-{1\over V}\int_X
\phi\o_0^n
=
J_{\o_0}(\phi)
-
{1\over V}\int_X\phi\o_0^N
\ee
For $\o_0\in c_1(X)$, the critical points of both the functional $F_{\o_0}(\phi)$
defined by
\be
F_{\o_0}(\phi)=F_{\o_0}^0(\phi)-\ln ({1\over V}\int_Xe^{h_{\o_0}-\phi}\o_0^n),
\ \ Ric(\o_0)-\o_0\equiv {\sqrt{-1}\over 2\pi}\p\bar\p h_{\o_0}
\ee
and the Mabuchi functional $\nu_{\o_0}(\phi)$ give K\"ahler-Einstein metrics.
Thus any direct relation between stability and the behavior of these functionals
would be valuable. Some progress in this direction can be found in \cite{tian97}\cite{zhang}\cite{tian00}\cite{PS02}.

\medskip

In the case of complex curves, the asymptotic behavior of $F_{\o_0}^0$ has been derived
by Paul \cite{P}. He obtained the formula (\ref{paul}) listed below.
This formula turns out to be related to
Mumford's algebraic formula \cite{mu}, and confirms in this case the close relation between energy functionals
and Chow-Mumford stability. The purpose of the
present paper is to determine, also in the case of
curves, the asymptotic
behavior of the Mabuchi energy functional $\nu_{\o_0}(\phi)$.

\medskip
Our set up is the following.
Let $L$ be a very ample line bundle on
$X$. Then $X$ can be imbedded into ${\bf CP}^N$ by
\be
X\ni z\to [S_0(z),\cdots,S_N(z)]
\ee
where $S_0(z),\cdots,S_N(z)$ is a basis for
the space $H^0(X,L)$ of holomorphic sections of $L$.
The line bundle $L$ corresponds then to the restriction
to $X$ of the hyperplane bundle ${\cal O}_{{\bf CP}^N}(1)$
over ${\bf CP}^N$.
Stability in the sense of geometric invariant theory
is a property of the action of $SL(N+1)$ acting on
$H^0(X,L)$ by change of bases $[S_0(z),\cdots, S_N(z)]$. It suffices actually to
consider one-parameter subgroups
$\sigma_t\subset SL(N+1,{\bf C})$, which can be assumed to act diagonally
\be
\sigma_t\cdot S=(t^{a_0}S_0(z),\cdots,t^{a_N}S_N(z)),
\ \
a_0+\cdots+a_N=0.
\ee
Under this action, $X$ acquires a one-parameter family of
K\"ahler metrics
\be
\o_t
={\sqrt{-1}\over 2\pi}\p\bar\p
\ln {||\sigma_t\cdot S||^2},
\ \
||\sigma_t\cdot S||^2=\sum_{j=0}^N|t|^{2a_j}|S_j(z)|^2,
\ee
which are the restrictions to
$\sigma_t\cdot X$ of the Fubini-Study metric on ${\bf CP}^N$. 
In terms of potentials,
$\o_t=
\o_0+{\sqrt{-1}\over 2\pi}\p\bar\p\,\phi$,
where $\o_0=\sqrt{-1}\p\bar\p\ln||S||^2$
and $\phi$ is defined by
\be
\phi(z)
=\ln {||\sigma_t\cdot S||^2\over||S||^2}
=\ln {\sum_{j=0}^N|t|^{2a_j}|S_j(z)|^2
\over
\sum_{j=0}^N|S_j(z)|^2},
\ee
All the expressions above
are valid in any local trivialization $z$ of the line bundle $L$,
with $|S_j(z)|$ being just the absolute value of the complex number
$S_j(z)$ in such a trivialization.
The issue which we address here is the
asymptotic behavior of $\nu_{\o_0}(\phi)$ as $|t|\to 0^+$.
Our main result is the following, with the notation to be explained in detail
in \S 2:

\bigskip

{\bf Theorem 1} {\it Assume that $X$ has dimension $n=1$. Then
the asymptotic behavior of the Mabuchi energy functional under
the action of a one-parameter subgroup of $SL(N+1)$ is given by}
\be
\nu_{\o_0}(\phi)
=
\,
\ln {1\over |t|}
\,
\sum_{zeroes\ of\ S_N}
{1\over V}\big\{2q_0-\mu\sum_{\al=1}^M
p_{\al}^2(m_{\al}-m_{\al+1})\big\}
+O(1)
\ee
{\it Here $q_0$, $p_{\al}$, and $m_{\al}$ refer to the Newton diagram
of $\s_t\cdot S$ at each fixed zero of $S_N(z)$.}

\bigskip
The main idea in the proof of Theorem 1 is a decomposition
of a neighborhood of the zeroes of $S_N(z)$ into annuli,
on each of which the integrand can be simplified to
essentially the Green's function for the Laplacian.
The method works well because the only singularity is a single
power of $\ln 1/|t|$. More general methods for the evaluation
of integrals depending on a holomorphic parameter $t$
are in \cite{PS00}. Related methods for estimating integral
operators rather than the scalar integrals discussed here 
can be found in \cite{ps,pss}.

\section{Degeneracies of $\sigma_t\cdot S(z)$ and of $\omega_t$}
\setcounter{equation}{0}

It is convenient to introduce the following notation
\be
\phi=\ln{|\s S|^2\over |S|^2}-2a_N\ln{1\over|t|},
\ \ \
|\s S|^2=
\sum_{j=0}^N|t|^{2q_j}|S_j(z)|^2,
\ \ \
|S|^2=
\sum_{j=0}^N|S_j(z)|^2,
\ee
with the exponents $q_j$ given by
\be
q_j=a_j-a_N\geq 0.
\ee
The asymptotics of the energy functionals will be given by
the degeneracies of $S_j(z)$ near the zeroes of $S_N(z)$. To describe
them, let $z_0\in X$ be a zero
of $S_N(z)$, with $z_0=0$ in the local
trivialization $z$ of $X$.
Let
\be
S_j(z)=u_jz^{p_j}+O(z^{p_j+1}),
\ \ u_j\not=0
\ee
near $0$. The Newton diagram
of $\s S$ is defined to be the Newton diagram in ${\bf R}_+\times {\bf R}_+$
defined by the points $(p_j,q_j)$.
Recall that the Newton diagram defined by a set of points $(p_j,q_j)$
is the convex hull of the union of all upper quadrants
$\{(p,q);\ p\geq p_j\ {\rm and}\ q\geq q_j\}$
with corner at $(p_j,q_j)$.
Let $\{V_{\alpha}\}_{\al=0}^M$ be the set of vertices of the Newton diagram.
It is possible for several indices $j$ to produce the
same vertex $(q_j,p_j)=V_{\al}$, in which case we also introduce
the notation $|u_{\alpha}|^2=\sum_{(p_j,q_j)=V_{\al}}|u_j|^2$.
Let $\{m_{\al}\}_{\al=1}^M$ be the set
of slopes, listed in {\it decreasing} order
\bea
&&
q_{\al-1}-q_{\al}+m_{\alpha}(p_{\al-1}-p_{\al})=0,
\ \ V_{\al}=(p_{\al},q_{\al}),
\nonumber\\
&&
m_1>m_2\cdots>m_M.
\eea
We observe that there is a vertex $V_0$ on the $q$-axis, since
$||S(z)||^2$ is never $0$. Similarly, by construction of the $q_j$,
there is a vertex $V_M$ on the $p$-axis. This implies in particular
that all the slopes $\{m_{\alpha}\}_{\al=1}^M$ are finite and strictly positive.
It is convenient to introduce also
the trivial slopes $m_0=+\infty$ and $m_{M+1}=0$ of the Newton diagram.

\medskip

We also need the degeneracies of $\o_t$.
In a local trivialization, we may write
\be
\label{ot}
\o_t=
\p_z\p_{\bar z}
\ln |\s S|^2{idz\wedge d\bar z\over 2\pi}
=
{\sum_{0\leq j<k\leq N}|t|^{2(q_j+q_k)}
|S_k\p_zS_j-S_j\p_zS_k|^2
\over
2(\sum_{j=0}^M|t|^{2q_j}|S_j(z)|^2)^2}
{idz\wedge d\bar z\over 2\pi}
\ee

\medskip
The basic idea in our approach is to decompose a fixed
neighborhood $\{|z|<1\}$ of $0$ into annuli of the
form
\be
\label{annulus}
|t|^{m_{\al}}\leq |z|<|t|^{m_{\al+1}},
\ \ 0\leq\al\leq M.
\ee
The advantage is that the term $|t|^{2q_\al}|z|^{2p_\al}$ dominates all others in this annulus. More precisely, we have
\be
\label{size}
|t|^{2q_j}|z|^{2p_j}\leq |t|^{2q_{\al}}|z|^{2p_{\al}},
\ \ 0\leq j\leq M,
\ee
with the terms $|t|^{2q_j}|z|^{2p_j}$
matching for $|z|=|t|^{m_\al}$ and $(p_j,q_j)$ on the face $F_{\al}^-$ linking
$V_{\al-1}$ to $V_{\al}$, and for $|z|=|t|^{m_{\al+1}}$
for $(p_j,q_j)$ on the face
$F_{\al}^+$ linking $V_{\al}$ to $V_{\al+1}$:
\bea
\label{relative}
|t|^{2q_j}|z|^{2p_j}
&=&
|t|^{2q_{\al}}|z|^{2p_{\al}},
\ \ {\rm when}\ |z|=|t|^{m_{\al+1}}, \ \ (p_j,q_j)\in F_{\al}^+
\nonumber\\
|t|^{2q_j}|z|^{2p_j}
&=&
|t|^{2q_{\al}}|z|^{2p_{\al}},
\ \ {\rm when}\ |z|=|t|^{m_{\al+1}}, \ \ (p_j,q_j)\in F_{\al}^-.
\eea
In particular, in the annulus (\ref{annulus}), we have
$c|t|^{2q_{\al}}|z|^{2p_{\al}}\leq
\sum_{j=0}^M|t|^{2q_j}|S_j(z)|^2
\leq C|t|^{2q_{\al}}|z|^{2p_{\al}}$
for suitable constants $c,C>0$, and
\bea
\label{approx}
\ln |\s S(z)|^2
&=&
\ln\,(|t|^{2q_{\al}}|z|^{2p_{\al}}|u_{\al}|^2)
+
O({\sum_{(p_j,q_j)\not=V_{\al}}|t|^{2q_j}|S_j(z)|^2
\over
|t|^{2q_{\al}}|z|^{2p_{\al}}})
\nonumber\\
{1\over |\s S|^2}
&=&
{1\over |t|^{2q_{\al}}|z|^{2p_{\al}}|u_{\al}|^2}+O({\sum_{k\not\in V_{\al}}|t|^{2q_k}|S_k|^2
\over (|t|^{2q_{\al}}|z|^{2p_{\al}})^2})
\eea

The following elementary lemma is useful in dealing with the type of error terms
arising in the preceding expansions:

\bigskip
{\bf Lemma 1}. {\it We have the following asymptotics
\bea
&&
\int_{|t|^{m_{\al}}\leq |z|<
|t|^{m_{\al+1}}}
|t|^{2(q_1+\cdots+q_M-Mq_{\al})}
|z|^{2(p_1+\cdots+p_M-Mp_{\al})}{d^2z\over 2\pi |z|^2}
\nonumber\\
&&
\quad\quad
=
\delta_{V_{\alpha}}(j_1,\cdots,j_M)
\ln\,{1\over |t|}\,
(m_{\al}-m_{\al+1})+O(1)
\eea
where the Dirac function $\delta_{V_{\alpha}}(j_1,\cdots,j_M)$
is defined to be $1$ if all the points $(p_j,q_j)$, $j=1,\cdots, M$, coincide with the vertex $V_{\al}$ and $0$ otherwise.
A similar statement holds for the interval of integration
$|t|^{m_{\al-1}}\leq |z|<|t|^{m_{\al}}$, with
$m_{\al}-m_{\al+1}$ replaced by $m_{\al-1}-m_{\al}$
on the right hand side.}

\bigskip
{\it Proof of Lemma 1}. If $p_1+\cdots+p_M-Mp_{\al}=0$,
then $q_1+\cdots+q_M-Mq_{\al}\geq 0$ since all points $(p_j,q_j)$
are in the Newton diagram. Since $q_j-q_{\al}
+m_{\al}(p_j-p_{\al})\geq 0$ for each $j$,
we can have equality only
if all points $(p_j,q_j)$ are on the faces of the Newton diagram.
But if they are on the faces, then $q_j-q_{\al}=0$ since all
points on a face must be on one side of a vertex.
We can now evaluate the integral of $d^2z/2\pi|z|^2$ and
verify the desired formula in this case. When
$p_1+\cdots+p_M-Mp_{\al}\not=0$, the integral can be evaluated
directly, giving
\be
|t|^{2\sum_{j=1}^M(q_j-q_\al)+m_{\alpha+1}(p_j-p_\al)}
-
|t|^{2\sum_{j=1}^M(q_j-q_\al)+m_{\alpha}(p_j-p_\al)}
\ee
up to a multiplicative constant. This is $O(1)$. Q.E.D.

\section{Asymptotics for the Futaki functional}
\setcounter{equation}{0}

The asymptotics of the Mabuchi functional will be derived by
combining the asymptotics of the various terms in its
definition (\ref{mabuchi}). A first term is proportional to
$F_{\o_0}^0$, the asymptotics of which have been
derived by S. Paul in \cite{P}, as we noted earlier.
The formula obtained in \cite{P} is the following
\be
\label{paul}
F_{\o_0}^0(\phi)
=
\big\{2a_N+
{1\over V}\sum_{\al=1}^M
p_{\al}^2(m_{\al}-m_{\al+1})\big\}
\,
\ln{1\over |t|}
+O(1)
\ee
We give a different derivation of this result
now, for the convenience of the reader
and to illustrate the decomposition (\ref{annulus}).
The easiest way is to use the expression (\ref{aubinyau})
of $F_{\o_0}^0$
in terms of the $J_{\o_0}$
functional. First, we observe that
\be
{1\over V}
\int_X\phi\,\o_0
=
-2a_N\ln{1\over |t|}
-{1\over V}\int_X\ln{|S_0(z)|^2+\cdots+|S_N(z)|^2
\over
|t|^{2q_0}|S_0(z)|^2+\cdots+|S_N(z)|^2}\o_0
\ee
The integrand is between $0$ and
$\ln |S_0|^2+\cdots+|S_N(z)|^2/|S_N(z)|^2$.
This last expression is independent of $t$, and has only logarithmic singularities in $z$. Thus its integral is finite, and $O(1)$ as
$t\to 0$. It remains to estimate
\be
J_{\o_0}=
{\sqrt{-1}\over 4\pi V}
\int_X\p\phi\wedge\bar\p\phi
\ee
Evidently, we can restrict the integral to the
region $|z|\leq 1$ around a fixed zero of $S_N(z)$.
Since $|\p_z\ln |\s S|^2|\leq C|z|^{-1}$
and $|\p_z\ln |S|^2|$ is uniformly bounded, we may
replace $\p \phi\wedge\bar\p \phi$ by
$\p\ln |\s S|^2\wedge\bar\p\ln |\s S|^2$
in the integrand. Thus the desired formula
(\ref{paul}) is a consequence
of the following
\be
\int_{|t|^{m_{\al}}\leq |z|<|t|^{m_{\al+1}}}
\p_z\ln |\s S|^2\p_{\bar z}\ln |\s S|^2
{d^2z\over 2\pi}
=
\big(\ln{1\over |t|}\big)
\,p_{\alpha}^2(m_{\al}-m_{\al +1})
+O(1)
\ee
\bigskip
Now the term
\be
\p_z |\s S|^2
\p_{\bar z}|\s S|^2
=
\sum |t|^{2(q_j+q_k)}\p_zS_j\,S_j^*\p_{\bar z}S_k^*S_k
\ee
can be expanded near $0$ as
\be
p_{\al}^2|t|^{4q_{\al}}|z|^{4p_{\al}-2}(|u_{\al}|^2+O(|z|))
+
O\big(\sum_{(p_j,q_j;p_k,q_k)\not= (V_\al;V_\al)}
|t|^{2(q_j+q_k)}|z|^{2(p_j+p_k)-2}\big)
\ee
Combined with the expansion (\ref{size}) for $|\s S|^{-2}$,
this implies
\be
{\p_z|\s S|^2
\p_{\bar z}|\s S|^2
\over |\s S|^4}
=
{p_{\al}^2+O(|z|)\over |z|^{2}}
+
O({1\over |t|^{4q_\al}|z|^{4p_{\al}}}
\sum_{(p_j,q_j;p_k,q_k)\not= (V_\al;V_\al)}
|t|^{2(q_j+q_k)}|z|^{2(p_j+p_k)-2})
\ee
Integrating over the annulus $|z|^{m_{\al}}
\leq |z|<|t|^{m_{\al+1}}$, we find that
the above error terms all contribute $O(1)$ by Lemma 1, since at least
one of the points $(p_j,q_j)$ does not coincide with the vertex $V_{\al}$.
Again by Lemma 1, the main terms give the desired asymptotics.
Q.E.D.

\section{Asymptotics for the Mabuchi functional}
\setcounter{equation}{0}

The main additional difficulty in the Mabuchi functional is
the occurrence of the term $\p_z\p_{\bar z}\phi$,
which cannot be avoided as above in the case
of $F_{\o_0}^0$ for complex curves.
We shall handle such terms by Green's formula, combined with suitable
approximations valid in the annulus (\ref{annulus}).
As a warm up, we show how this technique works in the case of $F_{\o_0}^0$,
left under the form
\be
\label{F2}
F_{\o_0}^0(\phi)
=-{1\over 2V}\int_X\phi(\o_0+\o_\phi)
=-{\sqrt{-1}\over 4\pi V}\int_X\phi\,\p\bar\p\phi-{1\over V}\int_X\phi\o_0
\ee
The following simple lemma is useful:

\bigskip

{\bf Lemma 2.} {\it Let $T_1(z),\cdots,T_M(z)$ be holomorphic functions
in a neighborhood of $0$ which do not all vanish
identically. Then there exists a neighborhood ${\cal O}$
and a constant $C$ so that
\be
|\p_z\p_{\bar z}\,\ln\, (\sum_{l=1}^M\e^{2q_l}|T_l(z)|^2)|
\leq C{1\over |z|^2}
\ee
for all $z\in {\cal O}$, and all $\e>0$.}

\bigskip
{\it Proof of Lemma 2.} The left hand side can be bounded by
a linear combinations of terms of the form $\e^{2(q_j+q_k)}|T_j|^2
|\p_zT_k|^2(\sum_{j=1}^M\e^{2q_j}|T_j(z)|^2)^{-2}$. The desired estimate
follows from the estimate $|\p_zT_j(z)|\leq C_j|z|^{-1}|T_j(z)|$,
for constants $C_j$ independent of $\e$ and $z\in {\cal O}$. Q.E.D.

\subsection{The Futaki functional revisited}

We return to the asymptotics of the expression (\ref{F2}).
The term $\phi\,\o_0$ has been treated before, so we concentrate
on the term $\phi\p\bar\p\phi$.
It is easy to see that we can again localize to a neighborhood
$\{|z|<1\}$ of a zero of $S_N(z)$, and on that neighborhood,
replace $\phi(z)$ by the following more convenient
function $\psi(z)$ defined by
\be
\label{psi}
\psi(z)=\ln\,|\s S|^2
=
\ln\,(\sum_{j=0}^M|t|^{2q_j}|S_j(z)|^2)
\ee
Consider next the approximation of $\psi=\ln\,|\s S|^2$ given by (\ref{approx}).
In view of Lemmas 1 and 2, the error terms in that
approximation give rise only to bounded terms. Thus we can write
\be
\int_{|t|^{m_{\al}}\leq |z| <|t|^{m_{\al}}}
\psi\,
\p_z\p_{\bar z}\psi\,{d^2z\over 2\pi}
=
\int_{|t|^{m_{\al}}\leq |z| <|t|^{m_{\al}}}\ln (|t|^{2q_{\al}}|z|^{2p_{\al}}|u_{\al}|^2)
\p_z\p_{\bar z}\psi\,{d^2z\over 2\pi}
+O(1)
\ee
The Green's formula can now be applied to the integral on the right hand side
\bea
&&
\int_{|t|^{m_{\al}}\leq |z| <|t|^{m_{\al}}}\ln (|t|^{2q_{\al}}|z|^{2p_{\al}}|u_{\al}|^2)
(4\p_z\p_{\bar z}\psi)\,{d^2z\over 2\pi}
\nonumber\\
&&
\quad\quad
=
\oint_{|z|=|t|^{m_{\al+1}}}
\bigg\{\ln (|t|^{2q_{\al}}|z|^{2p_{\al}}|u_{\al}|^2)
{\p\over \p n}\psi
-
{\p\over \p n}\ln (|t|^{2q_{\al}}|z|^{2p_{\al}}|u_{\al}|^2)\psi
\big)\,\bigg\}{ds\over 2\pi}
\nonumber\\
&&
\quad\quad\quad
-\oint_{|z|=|t|^{m_{\al}}}
\bigg\{\ln
(|t|^{2q_{\al}}|z|^{2p_{\al}}|u_{\al}|^2)
{\p\over \p n}\psi\,ds
-
{\p\over \p n}\ln (|t|^{2q_{\al}}|z|^{2p_{\al}}|u_{\al}|^2)\psi
\big)\,\bigg\}{ds\over 2\pi}
\eea
Here we made use of the fact that $\p_z\p_{\bar z}\ln |z|^2=0$ on the annulus.
Now the asymptotics of all the terms in the above expression
at the boundaries $|z|=|t|^{m_{\al}+1}$ and $|z|=|t|^{m_{\al}}$ are very simple.
For $|z|=|t|^{m_{\al+1}}$, they are given respectively by
\bea
\psi(z)&=&
-2(q_{\al}+m_{\al+1}p_{\al})\ln {1\over |t|} +O(1)\nonumber\\
|z|{\p\over\p n}\psi(z)
&=&
{\sum_{l\in F_{\al}^+}2p_l|u_l|^2
\over
\sum_{l\in F_{\al}^+}|u_l|^2}+o(1)
\eea
and
\bea
\ln (|t|^{2q_{\al}}|z|^{2p_{\al}}|u_{\al}|^2)
&=&
-2(q_{\al}+m_{\al+1}p_{\al})\ln {1\over |t|} +O(1)
\nonumber\\
|z|{\p\over \p n}\ln (|t|^{2q_{\al}}|z|^{2p_{\al}}|u_{\al}|^2)
&=&
2p_{\al}
\eea
with similar asymptotics near $|z|=|t|^{m_{\al}}$.
Thus the contribution of the contour integrals involving ${\p\over\p n}\psi$
is
\be
-\ln{1\over |t|}\,\bigg[2(q_{\al}+m_{\al+1}p_{\al})
{\sum_{l\in F_{\al}^+}2p_l|u_l|^2
\over
\sum_{l\in F_{\al}^+}|u_l|^2}
-
2(q_{\al}+m_{\al}p_{\al})
{\sum_{l\in F_{\al}^-}2p_l|u_l|^2
\over
\sum_{l\in F_{\al}^-}|u_l|^2}\bigg]
\ee
Upon summation over all vertices $\al$, these contributions telescope to
$0$. Indeed, the right face $F_{\al-1}^+$ is the same as the left face
$F_{\al}^-$, and $q_{\al}+m_{\al}p_{\al}=q_{\al-1}+m_{\al}p_{\al-1}$,
since the vertex $(p_{\al-1},p_{\al-1})$ is on the face $F_\al^-$.
There are also no end terms in the telescoping series,
since there is no vertex to the left of $\al=0$ and to the right of $\al=M$.
We observe that this cancellation mechanism depends only on the fact
that $\ln (|t|^{2q_\al}|z|^{2p_\al}|u_{\al}|^2)$ and
$\ln (|t|^{2q_{\al-1}}|z|^{2p_{\al-1}}|u_{\al-1}|^2)$ have the same asymptotics
on $|z|=|t|^{m_{\al}}$. It does not depend on the exact value of ${\p \psi\over\p n}$,
and reflects the fact that the asymptotics of $F_{\o_0}^0$ are independent of the decomposition into annuli.

\medskip
It remains to consider the contributions of the contour integrals involving
$\psi$. Clearly they are for each $\al$
\be
\ln{1\over |t|}\,\bigg[
2p_{\al}2(q_{\al}+m_{\al+1}p_{\al})
-
2p_{\al}2(q_{\al}+m_{\al}p_{\al})\bigg]
=
-4p_{\al}^2(m_{\al}-m_{\al+1})
\,
\ln\,{1\over |t|}
\ee
agreeing with the earlier method.

\subsection{Proof of Theorem 1}

We treat now the Mabuchi functional $\nu_{\o_0}(\phi)$, concentrating
first on the new term
\be
\int_X\ln ({\o_t\over\o_0})\o_t
\ee
It still suffices to consider the
integrals over small neighborhoods of the isolated
zeroes of $S_N(z)$. This is because
$|\ln\,({\o_t\over\o_0})\,\o_t|\leq C$
if $S_N(z)$ is bounded away from $0$.
Thus, we may restrict to a neighborhood
$\{|z|<1\}$ of a zero of $S_N(z)$, and express the above
integral as
\be
\int_X\ln ({\o_t\over\o_0})\o_t
=
\int_X\ln ({\o_t\over\o_0})\p_z\p_{\bar z}\psi\,{2d^2z\over 2\pi},
\ee
in view of (\ref{ot}) and the definition (\ref{psi}) of $\psi(z)$.
We decompose again $\{|z|<1\}$ into
annuli $|t|^{m_{\al}}\leq |z|<|t|^{m_{\al+1}}$.
We consider separately three regions, when $|z|<|t|^{m_1}$,
when $|t|^{m_\al}\leq |z| <|t|^{m_{\al+1}}$ for
$1\leq \al\leq M-1$, and when $|t|^{m_M}\leq |z|<1$.

\bigskip

{\it The Region $|t|^{m_\al}\leq |z| <|t|^{m_{\al+1}}$}

\medskip

Let $\al$ be fixed with $1\leq \al\leq M-1$.
A first application of the Green's formula gives
\bea
&&
\int_{|t|^{m_{\al}}\leq |z|<|t|^{m_{\al+1}}}
(\ln\,{\o_t\over\o_0})\,
(4\p_z\p_{\bar z}\psi)\,{d^2z\over 2\pi}
-
\int_{|t|^{m_{\al}}\leq |z|<|t|^{m_{\al+1}}}
(4\p_z\p_{\bar z}\ln\,{\o_t\over\o_0})\,
\psi\,{d^2z\over 2\pi}
\nonumber
\\
&&
=
\oint_{|z|=|t|^{m_{\al+1}}}(\ln\,{\o_t\over\o_0})\,{\p\over\p n}\psi\,
{ds\over 2\pi}
-
\oint_{|z|=|t|^{m_{\al+1}}}({\p\over\p n}\ln\,{\o_t\over\o_0})
\,
\psi\,{ds\over 2\pi}
\nonumber
\\
&&
-
\oint_{|z|=|t|^{m_{\al}}}(\ln\,{\o_t\over\o_0})\,{\p\over\p n}\psi\,
{ds\over 2\pi}
+
\oint_{|z|=|t|^{m_{\al}}}({\p\over\p n}\ln\,{\o_t\over\o_0})
\,
\psi\,{ds\over 2\pi}
\eea

The second double integral on the left hand side can be simplified,
using the approximation for $\psi(z)$
in (\ref{approx}). Indeed, in local coordinates,
$\o_t$ is given by the expression (\ref{ot}). Lemma 2 implies then
\be
|\p_z\p_{\bar z}\,\ln\,
({\o_t\over\o_0})|
\leq C(1+{1\over |z|^2})
\ee
for some constant independent of $z$ and $|t|$.
Applying Lemma 1, we see that the error terms in the approximation
(\ref{approx}) for $\psi$ contribute only $O(1)$ terms to
the integral
\be
\int_{|t|^{m_{\al}}\leq |z|<|t|^{m_{\al+1}}}
(4\p_z\p_{\bar z}\ln\,{\o_t\over\o_0})\,
\psi\,{d^2z\over 2\pi}
=
\int_{|t|^{m_{\al}}\leq |z|<|t|^{m_{\al+1}}}
(4\p_z\p_{\bar z}\ln\,{\o_t\over\o_0})
\,\ln(|t|^{2q_{\al}}|z|^{2p_{\al}}|u_{\al}|^2)
{d^2z\over 2\pi},
\ee
up to $O(1)$ terms. Now apply Green's theorem again to
this new integral
\bea
&&
\int_{|t|^{m_{\al}}\leq |z|<|t|^{m_{\al+1}}}
(4\p_z\p_{\bar z}\ln\,{\o_t\over\o_0})
\,\ln(|t|^{2q_{\al}}|z|^{2p_{\al}}|u_{\al}|^2)
{d^2z\over 2\pi}
\\
&&
\quad
=
\oint_{|z|=|t|^{m_{\al+1}}}
\bigg\{({\p\over \p n}\ln\,{\o_t\over\o_0})
\ln (|t|^{2q_{\al}}|z|^{2p_{\al}}|u_{\al}|^2)
-
(\ln\,{\o_t\over\o_0})\,
{\p\over\p n}\ln (|t|^{2q_{\al}}|z|^{2p_{\al}}|u_{\al}|^2)\bigg\}
\,{ds\over 2\pi}
\nonumber\\
&&
\quad
-\oint_{|z|=|t|^{m_{\al}}}
\bigg\{({\p\over \p n}\ln\,{\o_t\over\o_0})
\ln (|t|^{2q_{\al}}|z|^{2p_{\al}}|u_{\al}|^2)
-
(\ln\,{\o_t\over\o_0})\,
{\p\over\p n}\ln (\e^{2q_{\al}}|z|^{2p_{\al}}|u_{\al}|^2)\bigg\}
{ds\over 2\pi}
\nonumber
\eea
Here we have exploited the fact that the function $\ln\,|z|^2$
is harmonic in the annulus. Altogether, we have then
\bea
\int_{|t|^{m_{\al}}\leq |z|<|t|^{m_{\al+1}}}
\ln\,{\o_t\over\o_0}\,(4\p_z\p_{\bar z}\psi)\,{d^2z\over 2\pi}
&=&
\oint_{|z|=|t|^{m_{\al+1}}}
\bigg\{(\ln{\o_t\over\o_0}){\p\over\p n}(\psi-
\ln(\e^{2q_{\al}}|z|^{2p_{\al}}|u_{\al}|^2))
\nonumber\\
&&
\quad
-
({\p\over\p n}\ln{\o_t\over\o_0})\,(\psi-
\ln(|t|^{2q_{\al}}|z|^{2p_{\al}}|u_{\al}|^2))\bigg\}{ds\over 2\pi}
\nonumber\\
&&
-\oint_{|z|=|t|^{m_{\al}}}
\bigg\{(\ln{\o_t\over\o_0}){\p\over\p n}(\psi-
\ln(\e^{2q_{\al}}|z|^{2p_{\al}}|u_{\al}|^2))
\nonumber\\
&&
\quad
-
({\p\over\p n}\ln{\o_t\over\o_0})\,(\psi-
\ln(|t|^{2q_{\al}}|z|^{2p_{\al}}|u_{\al}|^2))\bigg\}{ds\over 2\pi}
\nonumber
\eea
As before, the evaluation of the asymptotics of all the contour integrals
is easy. More specifically, we have near $|z|=|t|^{m_{al+1}}$
\bea
\psi-\ln(|t|^{2q_{\al}}|z|^{2p_{\al}}|u_{\al}|^2)
&=& O(1)
\nonumber\\
|z|{\p\over\p |z|}
(\psi-\ln(|t|^{2q_{\al}}|z|^{2p_{\al}}|u_{\al}|^2))
&=&
{\sum_{l\in F_{\al}^{+}}2p_l|u_l|^2
\over
\sum_{l\in F_{\al}^+}|u_l|^2}-2p_{\al}
\eea
Here we have denoted the left and right (closed)
faces of the Newton
diagram meeting at the vertex $(p_{\al},q_{\al})$
by $F_{\al}^-$ and $F_{\al}^+$
respectively.
Similarly, we have near $|z|=|t|^{m_{\al}}$
\bea
\psi-\ln(|t|^{2q_{\al}}|z|^{2p_{\al}}|u_{\al}|^2)
&=& O(1)
\nonumber\\
|z|{\p\over\p |z|}
(\psi-\ln(|t|^{2q_{\al}}|z|^{2p_{\al}}|u_{\al}|^2))
&=&
{\sum_{l\in F_{\al}^{-}}2p_l|u_l|^2
\over
\sum_{l\in F_{\al}^-}|u_l|^2}-2p_{\al}
\eea
Again, in these formulas, the exact value of $|z|{\p\over\p n}\psi$
will not be necessary because of its cancellation upon
summing in $\al$, but we listed it for the sake of completeness.
It is also easy to see that
\be
|z|{\p\over\p |z|}
\ln\,{\o_t\over\o_0}=O(1)
\ \ {\rm for}\ |z|=|t|^{m_{\al+1}}
\ or \
|z|=|t|^{m_{\al}}
\ee
In view of the above asymptotics
for $\psi-\ln(|t|^{2q_{\al}}|z|^{2p_{\al}}|u_{\al}|^2)$,
it is clear that the terms involving
$|z|{\p\over\p |z|}
(\ln\,{\o_t\over\o_0})$ only contribute $O(1)$. Thus it remains only
to derive the asymptotics for $\ln {\o_t\over\o_0}$ at $|z|=|t|^{m_{\al+1}}$
\be
\ln\,{\o_t\over\o_0}
=-\lambda_{\al+1}\,\ln{1\over |t|}+O(1)
\ee
To determine $\lambda_{\al}$, we recall the relative size
of $|t|^{2q_j}|z|^{2p_j}$ at $|z|=|t|^{m_{\al}}$ and
$|z|=|t|^{m_{\al+1}}$ given in (\ref{relative}).
It follows that the dominating terms in the expression
(\ref{ot}) for $\p_z\p_{\bar z}\ln |\s S|^2$
correspond respectively to $(j,k)=(\al,\al+1)$
and $(j,k)=(\al-1,\al)$ for
$|z|=|t|^{m_{\al+1}}$ and $|z|=|t|^{m_{\al}}$.
This implies
\be
\lambda_{\al+1}
=-2\,m_{\al+1}
\ee
Thus we obtain
\bea
&&
\int_{|t|^{m_{\al}}\leq |z|<|t|^{m_{\al+1}}}
\ln\,{\o_t\over\o_0}
(\p_z\p_{\bar z}\psi){d^2z\over 2\pi}
\nonumber\\
&&
\quad
=
\ln{1\over |t|}
\bigg\{({\sum_{l\in  F_{\al}^{+}}p_l|u_l|^2
\over
\sum_{l\in F_{\al}^+}|u_l|^2}-p_{\al})m_{\al+1}
-
({\sum_{l\in F_{\al}^{-}}p_l|u_l|^2
\over
\sum_{l\in F_{\al}^-}|u_l|^2}-p_{\al})m_{\al}\bigg\}
+O(1)
\eea
This completes the analysis of the region
$|t|^{m_{\al}}\leq |z|<|t|^{m_{\al+1}}$, for
$1\leq\al\leq M-1$.

\bigskip
{\it The Region $0\leq |z|<|t|^{m_1}$}

\medskip
This region required no separate treatment in the earlier case of $F_{\o}^0$. In the present case, there would be a Dirac contribution
at the origin if
the function $\o_t/\o_0$
also vanishes there. However, this cannot happen since $z\to
[S_0(z),\cdots,S_N(z)]$ has been assumed to be an imbedding.
This implies that $\o_t$ cannot vanish, since it is the restriction
to a smooth subvariety of the Fubini-Study metric on ${\bf CP}^N$.

\bigskip
{\it The Region $|t|^{m_M}\leq |z|<1$}

\medskip
This region also requires a separate argument, since at its outer boundary
$|z|=1$, the absolute value $|z|$ is merely small, without being of the order
of a {\it positive} power of $|t|$.
We shall show that, as in the model case of $F_{\o}^0$,
the outer boundary $|z|=1$ does not contribute, i.e.
\be
\int_{|t|^{m_M}\leq |z|< 1}
\ln\,{\o_t\over\o_0}\,(\p_z\p_{\bar z}\psi)\,{d^2z\over 2\pi}
=
-\ln{1\over |t|}
\cdot ({\sum_{l\in  F_M^-}p_l|u_l|^2
\over
\sum_{l\in F_M^-}|u_l|^2}-p_M)m_M
+O(1)
\ee
To do so, recall that the term $|S_N(z)|^2$ occurs without
factors of $|t|$ in $\psi(z)$, i.e., $q_N=0$.
We consider two separate cases. In the first case,
there are other indices $j$, besides $j=N$,
for which $q_j=0$.
For two such $j<k$, the term $|S_j\p_zS_k-S_k\p_zS_j|^2$
occurs in the expression (\ref{ot}) for $\o_t$
without any factor of $|t|$.
We may also assume that it does not vanish near $|z|=1$.
Thus $\o_t$ is bounded away from $0$, uniformly in
$t$, for $|t|$ small enough. Since $\sum_{j=0}^N|t|^{2q_j}|S_j(z)|^2$
is also bounded away from $0$ for $|z|=1$, the function $\ln(\o_t/\o_0)$ is then
smooth and bounded uniformly. It follows that in the boundary terms
resulting from the double application of the Green's formula,
the boundary $|z|=1$ contributes only $O(1)$ terms.

\bigskip
In the other case, we may assume that $q_j>0$ for all $j< N$.
Without loss of generality, we may assume that the local coordinate $z$
for $X$ is chosen so that $S_N(z)=z^{p_N}$.
Since we have then
\be
\psi(z)-\ln |z|^{2p_N}
=
\ln(1+{\sum_{j=0}^{N-1}|t|^{2q_j}|S_j(z)|^2
\over |z|^{2p_N}})
\ee
it follows immediately that this expression and its derivatives
are smooth in $z$ for
$|z|$ near $1$, and bounded by a strictly positive power of $|t|$.
Since we also clearly have $|\ln (\o_t/\o_0)|\leq C\ln (1/|t|)$
for $|z|$ near $1$,
the contribution of the boundary integral $|z|=1$ is again
$O(1)$ in this case. This completes the estimates
for the region $|t|^{m_N}\leq |z|<1$.

\bigskip

The summation over all annuli $|t|^{m_{\al}}\leq |z|
<|t|^{m_{\al+1}}$ produces as before a telescoping sum.
All the contributions of $|z|{\p\phi\over\p n}$ cancel as before.
Thus we are left only with the terms
\be
\int_{|z|<1 }\ln{\o_t\over\o_0}\,
\o_t=
2\ln{1\over |t|}
\sum_{\alpha=1}^M
p_{\al}(m_{\al}-m_{\al+1})
+O(1)
\ee
The sum on the right hand side can be rearranged as
\be
\sum_{\al=1}^Mp_{\al}(m_{\al}-m_{\al+1})
=
\sum_{\al=1}^Mm_{\al}(p_{\al}-p_{\al-1})
=
\sum_{\al=1}^M(q_{\al-1}-q_{\al})=q_0
\ee
in view of the fact that the $m_{\al}$'s are the slopes of the
Newton diagram.
We can now complete the proof of Theorem 1.
The asymptotics of the remaining terms in the Mabuchi energy functional
are given by
\be
-{1\over V}\int_X\phi\,Ric(\o_0)
=2a_N\,(\ln{1\over|t|})\,{1\over V}\int_X Ric(\o_0)+O(1)
=2a_N\mu\,\ln{1\over|t|}+O(1)
\ee
and
\be
-\mu F_{\o_0}^0(\phi)
=
-\mu\big\{2a_N+{1\over V}
\sum_{\al=1}^Mp_{\al}^2(m_\al-m_{\al+1})\big\}
\,\ln{1\over |t|}+O(1).
\ee
Assembling all the terms gives the asymptotics
stated in Theorem 1. Q.E.D.

\section{Remarks}
\setcounter{equation}{0}

In \cite{PS02}, the Mabuchi functional
$\nu_{\o_0}(\phi)$ has been related to two suitable norms
$||\s_t\cdot Chow(X)||$ and $||\sigma_t\cdot Chow(X)||_{\#}$ for the
Chow vector of the variety $X$, together with a current term
associated to the singular locus $Z_s$ of the Chow variety.
This current term is of great interest, since
it encodes delicate geometric properties of the imbedding
of $X$ into ${\bf CP}^N$.
It is however difficult to evaluate directly. We wish to point out
that Theorem 1 gives an upper bound for the current term
in the case of complex curves $X$, since
it is known that the norm $||\s_t\cdot Chow(X)||$
is proportional to $F_{\o_0}^0(\phi)$,
and $\ln ||\cdot||_{\#}\leq \ln ||\cdot||+C$ for a suitable
constant $C$.

\vfill\break

\end{document}